\documentclass[3p]{elsarticle}

\usepackage{ulem}
\usepackage{color} 
\usepackage{psfrag} 
\usepackage{subfig}
\usepackage{hyperref}
\usepackage{csvsimple}
\usepackage{amsmath}
\usepackage{amsthm}
\usepackage{graphics}
 \usepackage{pstool} 
\usepackage%[disable]
{todonotes}

\usepackage{pgfplots}
\usepackage{pgfplotstable}
\usepackage{filecontents}
\usepackage{pgf}
\usepgfplotslibrary{colorbrewer}

\usepackage{tikz}
\usetikzlibrary{decorations.pathmorphing}
\usepackage{tikzscale}
\usepackage{tikz-3dplot}
\usepackage{graphicx}

\usetikzlibrary{3d}
\usetikzlibrary{spy}
\usetikzlibrary{snakes,arrows,shapes}    
\usetikzlibrary{spy}
\usetikzlibrary{backgrounds}  
\usetikzlibrary{decorations}
\usetikzlibrary{positioning}
\usetikzlibrary{patterns}
\usetikzlibrary{calc} 
\usetikzlibrary{external} % Is activated with the command '\tikzexternalize'.
\usepackage{booktabs}   
\usepackage{makecell} 
\usepackage{colortbl}   
\usepackage{xspace}
\usepackage{color}     
\usepackage{setspace}   
 
\usepackage{soul}
\usepackage{ulem}
\usepackage{currfile}
\usepackage{import}

\usepackage{floatrow}
\usepackage{multicol, blindtext}
\newfloatcommand{capbtabbox}{table}[][\FBwidth]

\pgfplotsset{compat=1.7}

\newcommand{\findmax}[3]{
    \pgfplotstablesort[sort key={#2},sort cmp={float >}]{\sorted}{#1}%
    \pgfplotstablegetelem{0}{#2}\of{\sorted}%
    \let #3=\pgfplotsretval%
}

\pgfplotscreateplotcyclelist{myCycleListBlackAndWhite}
{%
  {black, mark=o},
  {black, mark=square},
  {black, mark=triangle},
  {black, mark=diamond}, 
  {black, mark=pentagon}, 
  {black, mark=*},
  {black, mark=square*},
  {black, mark=triangle*},
  {black, mark=diamond*}, 
  {black, mark=pentagon*}, 
}

\definecolor{darkgreen}{rgb}{0,0.4,0} 
\definecolor{darkbrown}{rgb}{0.5, 0.396, 0.09}

\definecolor{c1}{rgb}{0.0, 0.4196078431372549, 0.6431372549019608}
\definecolor{c2}{rgb}{1.0, 0.5019607843137255, 0.054901960784313725}
\definecolor{c3}{rgb}{0.6705882352941176, 0.6705882352941176,
0.6705882352941176} \definecolor{c}{rgb}{0.34901960784313724, 0.34901960784313724, 0.34901960784313724}
\definecolor{c4}{rgb}{0.37254901960784315, 0.6196078431372549,
0.8196078431372549} \definecolor{c}{rgb}{0.7843137254901961, 0.3215686274509804, 0.0}
\definecolor{c5}{rgb}{0.5372549019607843, 0.5372549019607843,
0.5372549019607843} \definecolor{c}{rgb}{0.6352941176470588, 0.7843137254901961, 0.9254901960784314}
\definecolor{c6}{rgb}{1.0, 0.7372549019607844, 0.4745098039215686}
\definecolor{c7}{rgb}{0.8117647058823529, 0.8117647058823529,
0.8117647058823529}

\pgfplotscreateplotcyclelist{myCycleListColor}
{%
  {blue, mark=o},
  {red, mark=square},
  {darkbrown, mark=triangle},
  {darkgreen, mark=diamond},
  {black, mark=pentagon},
  {blue, mark=*},
  {red, mark=square*},
  {darkbrown, mark=triangle*},
  {darkgreen, mark=diamond*},
  {black, mark=pentagon*},
}

\pgfplotscreateplotcyclelist{myCycleListColorNoMarkers}
{
  {black},
  {blue, dashed},
  {red, dotted},
  {darkbrown, loosely dashed},
  {darkgreen, loosely dotted},
}

\pgfplotsset{every axis/.append style= 
              {
                font=\small,
                mark size=2,
                line width = 0.1,
                legend style={font=\small, mark size=3, draw=none, fill=none},
                legend cell align=left,
                cycle list name=myCycleListColor,
              }
            }
            
\pgfplotstableset
{
  every odd row/.style={before row={\rowcolor[gray]{0.9}}},
  every head row/.style=
  {
    before row=
    {%
      \toprule
    },
    after row=\midrule
  },
  every last row/.style=
  {
    after row=\bottomrule
  }
}

\newif\ifdrawboundingbox

\usetikzlibrary{external} % Is activated with the command
\tikzset{external/system call={pdflatex \tikzexternalcheckshellescape
-halt-on-error -interaction=batchmode -jobname "\image" "\texsource"}} 
\pgfkeys
{ 
  /pgf/number format/.cd,
  fixed,
  set thousands separator={\,},
  1000 sep in fractionals,
}

\newcolumntype{C}[1]{>{\centering\arraybackslash}m{#1}}
\newcolumntype{R}[1]{>{\raggedright\arraybackslash}m{#1}}
\newcolumntype{L}[1]{>{\raggedleft\arraybackslash}m{#1}}

\newcommand{\delete}[1]{\xspace}

\drawboundingboxtrue
\drawboundingboxfalse 

\usepackage{amsmath,amsfonts}
\usepackage{etoolbox}
\ifdef{\R}{
	\renewcommand{\R}{\mathbb{R}}
}{
\newcommand{\R}{\mathbb{R}}
}
\ifdef{\N}{\renewcommand{\N}{\mathbb{N}}}{\newcommand{\N}{\mathbb{N}}}
\ifdef{\C}{\renewcommand{\C}{\mathbb{C}}}{\newcommand{\C}{\mathbb{C}}}
\ifdef{\Z}{\renewcommand{\Z}{\mathbb{Z}}}{\newcommand{\Z}{\mathbb{Z}}}

\newcommand{\maximum}[1]{\ensuremath{  \max\left\{ #1 \right\} } }

\newcommand{\lnorm}[1]{{\left\vert\kern-0.25ex\left\vert\kern-0.25ex\left\vert #1 
    \right\vert\kern-0.25ex\right\vert\kern-0.25ex\right\vert}}

\newcommand{\refeq}[1]{(\ref{#1})}

\newcommand{\meinlemma}{Lemma \stepcounter{smeinTheorems}\the\value{smeinTheorems}}
\newcommand{\meinsatz}{Satz \stepcounter{smeinTheorems}\the\value{smeinTheorems}}
\newcommand{\meinedefinition}{Def.\ \stepcounter{smeinTheorems}\the\value{smeinTheorems}}

\newcounter{meineAussage} 
\newcommand{\opnorm}[1]{{\left\vert\kern-0.25ex\left\vert\kern-0.25ex\left\vert #1 
    \right\vert\kern-0.25ex\right\vert\kern-0.25ex\right\vert}}

\usepackage{svg}
\usepackage{algorithm}% http://ctan.org/pkg/algorithms
\usepackage{algpseudocode}% http://ctan.org/pkg/algorithmicx
\usepackage{gensymb}
\usepackage{todonotes}
\makeatletter
\renewcommand*{\ALG@name}{Procedure}
\makeatother
\usepackage[capitalise,noabbrev,sort&compress]{cleveref}
\usepackage{soul,xcolor}
\usepackage{gensymb}
\usepackage[T1]{fontenc}

\begin{document}

\begin{frontmatter}

\title{Two-level Method Part-scale Thermal Analysis of Laser Powder Bed Fusion Additive Manufacturing}

\author[PaviaAddress]{Massimo Carraturo}
%  \cortext[mycorrespondingauthor]{Corresponding author}
%  \ead{massimo.carraturo@unipv.it}
%
\author[GSSIAddress]{Alex Viguerie\corref{mycorrespondingauthor}}
  \cortext[mycorrespondingauthor]{Corresponding author}
     \ead{alexander.viguerie@gssi.it}

\author[PaviaAddress]{Alessandro Reali}
\author[PaviaAddress]{Ferdinando Auricchio}
%
%\address[PaviaAddress]{Department of Civil Engineering and Architecture, 
\address[GSSIAddress]{Gran Sasso Science Institute, 
viale F. Crispi 7, 67100 L`Aquila, Italy} 

\address[PaviaAddress]{Department of Civil Engineering and Architecture, 
University of Pavia,
via Ferrata 3, 27100 Pavia, Italy}

\newcommand{\publicationDate}{\today}

%\usepackage{fancyhdr}
%\date{}
%\fancypagestyle{plain}{%
%  \renewcommand{\headrulewidth}{0pt}%
%  \fancyhf{}%
%  \fancyfoot[L]{
%  \vspace{-1.25cm} 
%  \footnotesize{Preprint submitted to \textit{\journal{}}  \hfill
%  \publicationDate
%  \\
%  %The final publication is available at
%  % \url{http://myDOI} 
%  }} } 

% Abstract ---------------------------------------
\vspace{-1.5cm} 
\hrule 

\begin{abstract}
Numerical simulations of a complete laser powder bed fusion (LPBF) additive manufacturing (AM) process are extremely challenging or even impossible to achieve without a radical model reduction of the complex physical phenomena occurring during the process. However, even when we adopt reduced model with simplified physics, the complex geometries of parts usually produced by LPBF AM processes make this kind of analysis computationally expensive. 
In fact, small geometrical features - which might be generated when the part is design following the principal of the so-called design for AM, for instance, by means of topology optimization procedures - often require complex conformal meshes. Immersed boundary methods seem to offer a valid alternative to deal with this kind of complexity. The two-level method lies within this family of numerical methods and  presents a very flexible tool to deal with multi-scale problems. 
In this contribution, we apply the recently introduced two-level method to part-scale thermal analysis of LPBF manufactured components, first validating the proposed part-scale model with respect to experimental measurements from the literature and then applying the presented numerical framework to simulate a complete LPBF process of a topologically optimized structure, showing the capability of the method to easily deal with complex geometrical features.
\end{abstract}
 
%% Keywords ---------------------------------------
%\vspace{0.25cm}
%\noindent \textit{Keywords:} 
 
%\vspace{0.25cm}
%\hrule 

\end{frontmatter}

% Actual Content ---------------------------------
\sloppy

% % Appendix
%\appendix
\section{Introduction}\label{sec:introduction}
%% ======== General introduction  ===========

%Overview of LPBF process and their advantages ---------------------------------------

Laser powder bed fusion (LPBF) is an additive manufacturing (AM) process consisting of either a laser or an electron beam selectively melting a bed of metal powder in a layer-by-layer fashion. Such a technology has seen an exponential growth over the last decades thanks to its ability to produce components with small geometrical features covering a large range of scales \cite{king_laser_2015}. LPBF allows close-to-freeform production and has radically changed the approach to design, moving from the so-called design for manufacturing to functional design where the functionality of the component is optimized with very few constraints coming from the manufacturing process \cite{Gibson2015}. Several metal alloys can be processed by means of LPBF technology, the most common ones being iron, titanium, nickel, and aluminum alloys \cite{yap2015review}. In the present work, an LPBF process of stainless steel 316L is considered.

%Challenges in LPBF technology --------------------------------------------------------

Despite the above-mentioned potentiality, a widespread adoption of LPBF technology is hindered by the lack of standardization and repeatability associated with the complex physical nature of this manufacturing process. Nowadays, one of the main issues that the research community is trying to address is to understand the \textit{process-structure-property-performance} relationship occurring in LPBF products \cite{Zhang2020}. Several studies investigate the effects of LPBF process parameters at different scales \cite{Khairallah2016,Keller2017,Smith2016a,carraturo2020numerical}, from micro-structure to part-deflection, but a deterministic, holistic model has not yet been developed.

%Limitations of purely empirical process optimization ---------------------------------
Since the process is principally thermally-driven, most of the process-induced flaws in the final structures are influenced by the thermal history of the part, i.e., by the spatio-temporal evolution of the temperature field within the structure \cite{herzog2016additive,ghosh_single-track_2018}. In fact, a specific region of the structure can lead to higher residual heat compared to others and thus to locally varying micro-structural properties which directly affect the final performance (e.g., yield stress, fatigue life, etc.) of the component \cite{yang2020laser,benedetti2021architected}. Therefore, a purely empirical process parameter optimization might be ineffective and lead to very high costs.

%Advantages of accurate numerical thermal models --------------------------------------
Numerical modeling and simulations can play an important role enabling prediction of the thermal history in LPBF components without undergoing long and expensive \textit{trial-and-error} campaigns. Due to the extreme spatio-temporal ranges involved and the complexity of the physical process, LPBF simulations are potentially challenging from both a modeling and a numerical perspective. For detailed reviews on LPBF process simulations, interested readers are referred to~\cite{GALATI20181,YAN2018427,Schoinochoritis2017}. 

The Finite Element Method (FEM) is often used to perform numerical analysis of LPBF problems. Some authors have tried to employ adaptive mesh refinement and coarsening \cite{patil_novel_2014,denlinger_thermomechanical_2014,kollmannsberger_hierarchical_2018,carraturo2019a}; however, the most widespread solution is to adopt reduced/surrogate models of the thermal problem, which in the following we refer to as part-scale models.

Many different part-scale thermal models have been recently developed and validated, but this kind of research is far from being complete. These models play a key role in both thermo-mechanical and microstructural predictions. In fact, even if residual stresses and microstructural properties are generated at the local scale and have a strong dependency on the cooling rates, they are also influenced by the geometry of the structure due to the complex \textit{process-structure-property-performance} relationship in AM processes discussed above. 

%% ======== Literature review ===========

% Williams et al., ADDMA (2018)
The so-called pragmatic approach, first presented by \citet{WILLIAMS2018416} activates blocks of agglomerated layers at a given temperature (e.g., the melting temperature in \cite{WILLIAMS2018416} and a material-dependent relaxation temperature in \cite{yang_2019}). Results obtained using this method are in good agreement with experimental measurements of part deflection, but no validation is provided for the thermal model. Moreover, this method requires a conformal mesh. This may lead to very large systems when complex geometries are investigated, since the generated mesh may be over-refined in parts of the domain where critical geometrical features (e.g., sharp corners, small holes, etc.) are present.

% Peng et al., ADDMA (2018)
\citet{peng2018fast} develop a quasi-static thermo-mechanical model for the fast prediction of thermal distortion in LPBF processes using a thermal circuit network on a voxelized representation of the part domain. In such a model, Hoelzel and co-workers also neglect thermal inertia effects justifying this choice since the mechanical problem is mainly quasi-static at part-scale and no validation is provided for the thermal model.
 
% Neiva et al., FEA (2020)
A part-scale thermal model employing the virtual domain approximation is presented by \citet{NEIVA2020103343} allowing for a domain reduction, since heat losses through the powder and the base plate are modeled by means of an equivalent heat transfer coefficient. Such a model shows a relative error below $15\%$ with respect to temperature measurements obtained by means of thermocouples on an oblique square prism of $30\times 30\times 80$ mm$^3$ printed using Ti-6Al-4V.  

% Dugast et al., ADDMA (2020)
Another part-scale thermal model based on a highly-efficient matrix-free GPU computation is introduced by \citet{dugast2021part} leading to a speed-up factor higher than $\times$4 with respect to classical CPU-based implementations. The base plate temperature of a $20\times 20\times 20$ mm$^3$ cubic specimen is measured by means of a K-type thermocouple embedded in the center of the build plate. To and co-workers assume in this model adiabatic boundary conditions at the solid/powder interface and conduct a thorough sensitivity analysis with respect to the other calibrated parameters.

% Yavari et al., ADDMA (2021)
\citet{yavari2021thermal} propose a part-scale thermal model using graph-theory. Such a model is validated with respect to experimental thermal measurements for a 250 cm$^3$ volume impeller. A good agreement among the results obtained using the graph-theory approach and the experimental measurements is reported. However, the authors assume linear thermal parameters and their results strongly depend on the chosen number of nodes and
a gain factor which controls the rate of heat diffusion through the nodes. Both these parameters have been calibrated in a previous work \cite{yavari2020thermal} for the same material and AM machine as in \cite{yavari2021thermal}  but, since they have no direct physical meaning, turn out to be rather complex to be accurately calibrated.

% Carraturo et al., ADDMA (2020)
In \citet{carraturo2020modeling} an immersed boundary method, namely the Finite Cell Method, is employed to implement a part-scale thermo-mechanical model. The model is then validated with respect to experimental measurements on part deflection and residual strains carried out at the National Institute of Standards and Technology (NIST) for the AM 2018 Benchmark Series \cite{ambench2018} but no temperature measurements were provided to validate the accuracy of the predicted temperature field.

Avoiding a conformal mesh discretization by means of an implicit representation of the part geometry, the numerical method adopted in the present work follows the line of thought presented in \cite{CARRATURO2021102077,carraturo2020modeling}. However, in the present contribution, we introduce a \textit{local-global} approach, which may better capture the multi-scale feature of LPBF processes, employing different thermal models at the local and the global scale.  

%% ======== Novelty and objective of the present study  ===========

In particular, a part-scale thermal model is developed and experimentally validated using the two-level method first presented by \citet{viguerie2020fat} and analyzed in detail in \cite{VVBA2020}. Such a numerical approach is a local-global methodology which allows one to couple the results of a local, finer discretization  with a coarser, global representation of the problem solutions. One of the main advantages of this approach is its ability to keep both the local and the global discretization on a structured cartesian grid, making the mesh generation process computationally inexpensive and, at the same time, avoiding continuous re-meshing \cite{viguerie2021numerical}. In the two-level method, the scale separation can be achieved not only in space, but also in time. Moreover, the simple structure of both the local and the global mesh makes possible a full parallelization of the code. However, both these possibilities are left to future studies and not included in the present contribution.

\citet{williams2019situ} demonstrated that the dwell temperature, i.e., the temperature of the powder-bed upper surface at the instance of laser activation, plays a crucial role in the resulting melt-pool shape during the layer printing process, leading to different microstructure and porosity in the final part. 
Starting from this observation, we aim at developing a part-scale thermal model which - exploiting the multi-scale nature of the two-level method - delivers an accurate prediction of the dwell temperature. One of the key features of the proposed model is that the powder-bed is included in the computational domain. As demonstrated by the presented results, such an assumption allows us to easily calibrate the physical parameters of our model and, at the same time, to achieve a remarkable agreement with respect to the measured experimental data.     
 
The outline of this paper is as follows. In \cref{sec:governingEquations} the governing equations of the considered physical problem at the local and global scale are presented. \Cref{sec:numericalMethod} describes the numerical implementation details of the part-scale two-level approach. In \cref{sec:ResultsAndDiscussion}, the proposed methodology is first validated with respect to experimental measurements taken from \cite{williams2019situ} for an AM cylindrical component, and then is applied to simulate the complete LPBF process of a topologically optimized beam structure. Finally, in \cref{sec:conclusions} we draw the main conclusions and possible further extensions to the present contribution.

%% ======== Outline  ===========

%The outline of this paper is as follows. In~\cref{sec:governingEquations}, we present the set of governing equations used to describe the weakly-coupled thermo-mechanical problem.~\Cref{sec:experimentalSetup} briefly describes the setup of experimental measurements obtained at NIST. In~\cref{sec:numericalMethod}, we first shortly recall the main ideas underlying FCM and then detail its implementation to solve part-scale LPBF AM processes. In~\Cref{sec:resultsAndDiscussion},  we present and discuss our numerical results, comparing them with the experimental measurements, while we draw conclusions and provide an outlook in~\cref{sec:conclusions}.

\section{Governing equations}\label{sec:governingEquations}
In the present section, a thermal model for LPBF process simulation based on the two-level method at part-scale  is presented. More details on thermal models in welding and additive manufacturing can be found by the interested readers in~\citep{goldak_2005,lindgren_numerical_2006,gougeBook2017}.
\subsection{Thermal problem}\label{ssec:thermalProblem}
Considering a material obeying Fourier's law of heat conduction, the LPBF thermal process can be modeled by means of the heat transfer diffusion equation:
\begin{equation}
\rho c \dfrac{\partial T}{\partial t}-\nabla\cdot\left( \kappa \nabla T \right)=Q
\label{eq:thermalEquation}
\end{equation}
where $T$ indicates the temperature field, $\rho=\rho(T)$ the temperature-dependent density of the material, $c=c(T)$ the temperature-dependent specific heat capacity, $\kappa=\kappa(T)$ the temperature-dependent thermal conductivity, $t$ the time, and $Q$ an equivalent heat source.
The latent heat term associated to the material phase-change is neglected in the presented work, due to its negligible influence when macroscopic, part-scale thermal effects are investigated (see, e.g.,~\citet{chiumenti_2017}). Boundary conditions will be discussed in detail in \Cref{ssec:thermalBC}.

\subsection{The Two-level Method}\label{ssec:TwoLevelMethod}

The two-level method is based on a re-formulation of the heat transfer equation \refeq{eq:thermalEquation} as two coupled problems referred in the following as the \textit{local} and the \textit{global} problem.
Let us consider a global domain $\Omega_+$ and a local domain $\Omega_-\subset \Omega_+$, such that, following the notation introduced in \cite{viguerie2020fat}, $T^+$ and $T^-$ are the local and the global solution respectively. The thermal conductivity $\kappa$ is defined as:
\begin{equation}
	\kappa := 
	\begin{cases}
		\kappa_+ & \text{in } \Omega_+\setminus \Omega_- \\
		\kappa_- & \text{in } \Omega_-
	\end{cases}
\end{equation}
with $\rho_+$, $\rho_-$, $c_+$, and $c_-$ defined analogously. We note some or all of these parameters may depend on time, space, or temperature. This dependence is understood if not explicitly denoted.
\par Let $\gamma$ denote the \textit{interface} between local and global problems, across which information is exchanged. For $\eta \in H^{-1/2}\left(\gamma \right)$, $\eta \xi_{\gamma} \in H^{-1}\left(\Omega_+\right) = H^1(\Omega_+)'$ is defined as the linear functional such that:
  \begin{align}\label{deltaDef}
  	\int_{\Omega_+} \left(\eta \xi_{\gamma}\right)w = \int_{\gamma} \eta w \quad\quad \forall w \in H_{\Gamma_D}^1\left(\Omega_+\right),
  \end{align} 
  \noindent consistent with the notation used in \cite{maury2001fat, bertoluzza2011analysis, bertoluzza2005fat}.
\par In order to achieve a more efficient numerical solution, which is necessary to deal with the LPBF thermal problem at part-scale, the original formulations of the global and local problem as presented by \citet{viguerie2020fat, VVBA2020, viguerie2021numerical} are slightly changed such as:
\begin{align}
	c_+\rho_+\dfrac{\partial T^+}{\partial t} - \nabla \cdot \left( \kappa_+\nabla T^+ \right) &= Q +\left(\kappa_+ - \kappa_-\right)\frac{\partial T^-}{\partial \boldsymbol{n} } \xi_{\gamma}  \qquad \text{in } \Omega_+ \label{eq:globalProblem}\\
	c_-\rho_-\dfrac{\partial T^-}{\partial t} - \nabla \cdot \left( \kappa_-\nabla T^- \right) &= Q \qquad \text{in } \Omega_-, \label{eq:localProblem}
\end{align}
with boundary conditions defined as in a following section. The problems \eqref{eq:globalProblem} and \eqref{eq:localProblem} are then solved in an iterative manner until convergence is satisfied.
\par When compared to the two-level method as implemented in \cite{viguerie2020fat, viguerie2021numerical}, this formulation greatly simplifies and reduces the cost of the assembly and evaluation process. More specifically, this novel formulation reduces the number of necessary integral evaluations, no longer needs the storage and evaluation of the derivatives of $\kappa$, and avoids the necessity of evaluating volumetric integrals non-conformal to the mesh discretization, which may ultimately require the application of an adaptive integration scheme, such as the ones adopted in other immersed boundary approaches (e.g., the Finite Cell Method \cite{Duster2017Ency} or the Cut Finite Element Method \cite{Burman2015}).
\par As a consequence, one incurs a minor loss of consistency in the following sense: Let $T^* \in H^1(\Omega_+)$ denote the exact solution of \eqref{eq:thermalEquation}, $T^+ \in H^1 (\Omega_+)$ the exact solution of \eqref{eq:globalProblem}, and $T^- (\Omega_-)$ the exact solution of \eqref{eq:localProblem}. Then the following consistency conditions are satisfied:
\begin{align}\label{globalConsistency}
    \int_{\Omega_+ \setminus \Omega_- } \left(T^* - T^+\right)v &= 0 \qquad \text{for all } v \text{ in } H^1\left( \Omega_+\right),\\
    \label{localConsistency}
        \int_{\Omega_- } \left(T^* - T^-\right)q &= 0 \qquad \text{for all } q \text{ in } H^1\left( \Omega_-\right).
\end{align}
However, the more stringent condition on $T^+$:
\begin{align}\label{strongConsistency}
    \int_{\Omega_+ \cap \Omega_- } \left(T^* - T^+\right)v &= 0 \qquad \text{for all } v \text{ in } H^1\left( \Omega_+\right),
\end{align}
which holds for the formulation shown in \cite{viguerie2020fat, viguerie2021numerical, VVBA2020}, no longer holds. In practice, we find this error is small, provided that the jump across $\gamma$ is not overly large (following the theory for the related Fat-boundary method given in \cite{bertoluzza2011analysis}). Indeed, for many problems, the scale separation between the global and local mesh is sufficiently large such that the discretization error dominates over this region, and any additional error incurred by the loss of consistency is negligible in comparison.

\subsection{Local and global problem boundary conditions}\label{ssec:thermalBC}
\begin{figure*}[h!]
	\centering
	\subfloat[Local domain BCs.\label{subfig:localBCs}]
	{
		\includegraphics[width=0.45\textwidth]{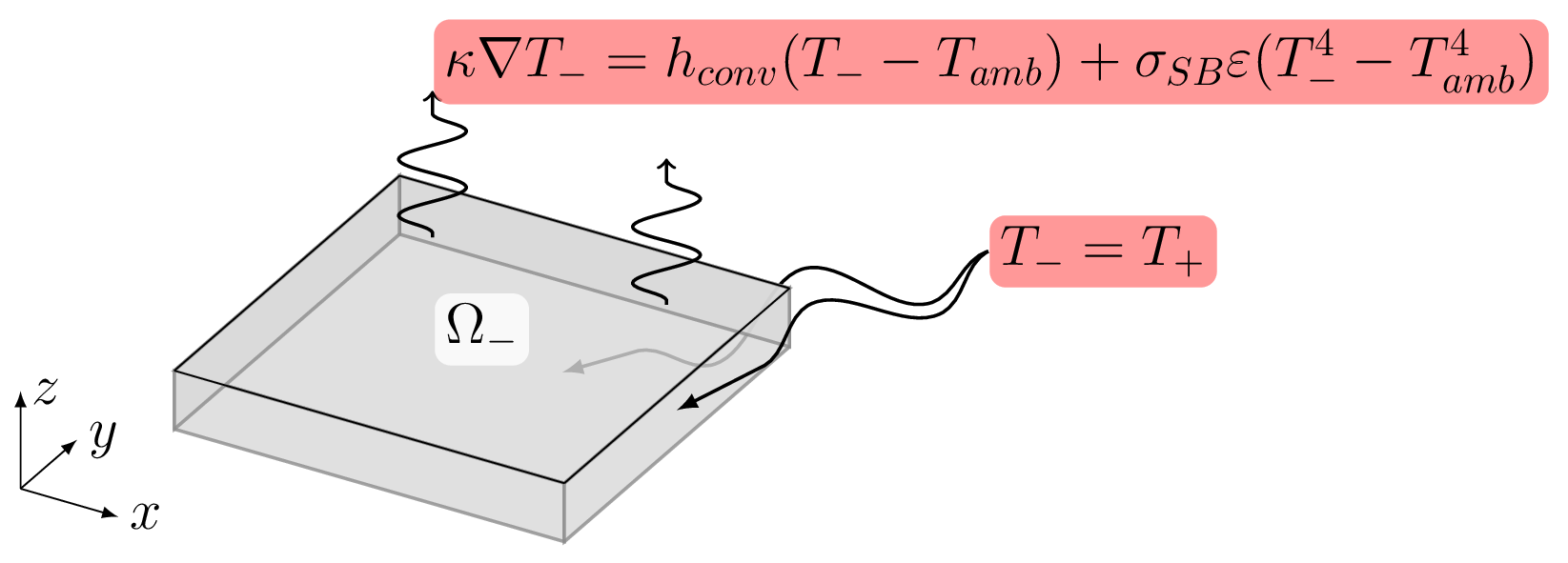}
	}
	\subfloat[Global domain BCs.\label{subfig:globalBCs}]
	{
		\includegraphics[width=0.45\textwidth]{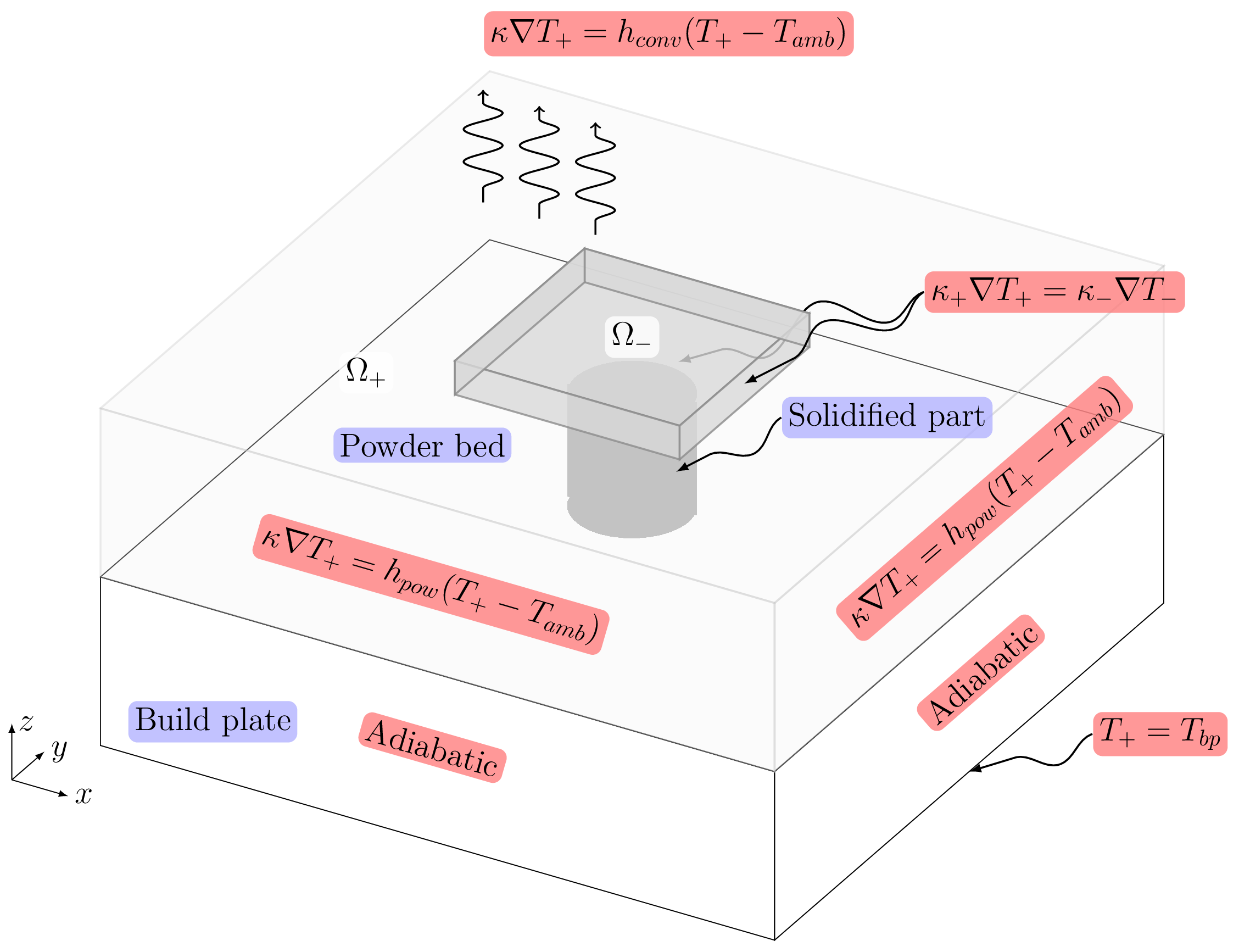}
	}    
	\caption{Thermal problem boundary conditions (BCs). In the local domain (a) the temperature continuity with the global solution is enforced by means of Dirichlet BCs on the lateral and lower surfaces whereas heat is dissipated through a convection and radiation heat flux on the upper surface. On the lower surface of the global domain (b) the temperature is fixed to the build plate temperature $T_{bp}$; adiabatic BCs are applied on the lateral surfaces of the build plate, while convection fluxes are imposed on the lateral and upper surfaces using different heat transfer coefficient values to model convection and conduction through the powder, respectively; finally, a zero-jump condition is imposed on the fluxes on the immersed local domain boundaries.  \label{fig:thermalBCs}}
\end{figure*}
To capture the multi-scale nature of the LPBF process, the two-level method distinguishes between a local domain $\Omega_-$ and a global domain $\Omega_+$, which we model differently in the proposed thermal model. \cref{fig:thermalBCs} depicts the boundary conditions (BCs) on the local and global domain. As can be observed, the local domain dissipates heat by convection and radiation through the upper surface (see \cref{subfig:localBCs}) by means of a heat loss flux term defined as follows:
\begin{equation}
	\kappa\nabla T_-= h_{conv}(T_- - T_{amb})+\sigma_{SB}\varepsilon (T_-^4 - T_{amb}^4),
\end{equation}
where $h_{conv}$ is the heat transfer coefficient by convection due to the inert gas (Argon) flow present in the chamber, $T_{amb}$ is the ambient temperature of the building chamber, $\sigma_{SB}=5.87\times 10^{-8}$ [W/m$^2$/K$^4$] is the Stefan-Boltzmann constant, and $\varepsilon$ the emissivity of the powder bed.
The global domain dissipates heat by conduction through the lateral surfaces of the powder bed and by convection through the upper surface. These two heat dissipation modes are modeled by means of two heat loss flux terms defined as:
\begin{align}
	\kappa\nabla T_+ &= h_{conv}(T_+ - T_{amb}),\\
	\kappa\nabla T_+ &= h_{pow}(T_+ - T_{amb}),
\end{align}
where $h_{pow}$ is the heat transfer coefficient by conduction through the powder. 

Most of these quantities are unknown and not readily measurable, and thus required a calibration procedure to be estimated.
Radiation effects are neglected in the global domain, since they play a minor role in regions far from the heat affected region, which we consider lying entirely within the local domain. As described in \cref{subfig:globalBCs}, global domain BCs are completed by adiabatic and Dirichlet BCs on the lateral and lower surfaces of the build plate, respectively. Finally, to avoid jumps among the local and the global domain temperature fields, a zero-jump temperature condition is imposed in the local problem on the lateral and bottom surfaces of $\Omega_-$, while, in the global problem, a zero-jump temperature flux condition is enforced on the lateral and bottom surfaces of $\Omega_-$ which are immersed in the global domain.

Due to the specific features of the two-level method as presented in \cite{viguerie2020fat}, the powder-bed is directly included in our thermal model. In general, immersed boundary methods do not require one to model the powder bed at part-scale \cite{carraturo2020modeling}. However, to impose temperature and flux continuity on the immersed local domain boundaries, the two-level method needs a physical solution on such boundaries, thus powder is included in our physical model. Modeling the powder surely increases the accuracy of the thermal model, but, on the other hand, it introduces potential numerical instabilities due to the discontinuous integration in the elements cut by the boundaries of the solidified domain.

\section{Numerical implementation}\label{sec:numericalMethod}
In the present work, we employ the two-level method \cite{viguerie2020fat,viguerie2021numerical}, an immersed boundary finite element method employing a local-global approach, together with an agglomerated layer activation procedure to solve the thermal problem defined in~\cref{sec:governingEquations} for an LPBF problem at part-scale. 

\subsection{Two-level algorithms}
The application of the two-level method may be done in several ways. Here, we consider two particular implementations, which we will refer to as the \textit{sequential two-level method} and as the \textit{parallel two-level method}, with details in the following. 
\par We consider an iteration $k$ at a generic time step $n$. We note that all variables indexed with $k$ refer to the current iteration, while those indexed with $n$ are from the preceding time step. For simplicity, we assume a standard Backward Euler time-stepping scheme, though we stress that the algorithm may be employed with any time integration technique. The two approaches are described in the following:

\noindent\textbf{Sequential two-level method:}
\newline \textbf{Step k:} Obtain local solution $T_k^-$:
\begin{align}\begin{split}\label{eq:seq_initialGuess}
		c_-|_{T_{k-1}^-} \rho_-|_{T_{k-1}^-} T_k^- -  \Delta t \nabla \cdot \left(\kappa_-|_{T_{k-1}^-} \nabla T_k^- \right) &= \Delta t Q + c_-|_{T_{k-1}^-} \rho_-|_{T_{k-1}^-} T_{n-1}^- \qquad \text { in } \Omega_-,  \\
		T_{k}^- &= T_{k-1}^+ \qquad \text{ on } \gamma.
\end{split}\end{align}
\newline \textbf{Step k+1:} Obtain global solution $T_k^+$:
\begin{align}\begin{split}\label{eq:seq_update}
		c_-|_{T_{k-1}^+} \rho_-|_{T_{k-1}^+} T_k^+ -  \Delta t \nabla \cdot \left(\kappa_-|_{T_{k-1}^+} \nabla T_k^+ \right) &= \Delta t Q + c_-|_{T_{k-1}^+} \rho_-|_{T_{k-1}^+} T_{n-1}^+ \\ &+ \Delta t \left(\kappa_+|_{T_{k-1}^+} - \kappa_- |_{T_k^-}\right)\frac{\partial T_k^-}{\partial \boldsymbol{n}} \xi _{\gamma}   \qquad \text { in } \Omega_+,  \\
		T_{k}^+ &= T_{0} \qquad \text{ on } \Gamma_D.
\end{split}\end{align}
\newline \textbf{Step k+2:} Relaxation: $T_k^+ = \theta T_{k}^+ + (1-\theta)T_{k-1}^+ $, $0<\theta \leq 1$.
\newline \textbf{Step k+3:} Check convergence. If reached, terminate. Else, return to step $k$.
\vspace{5mm}
\newline \textbf{Parallel two-level method:}
\newline \textbf{Step k:} Obtain local solution $T_k^-$:
\begin{align}\begin{split}\label{eq:parallel_initialGuess}
		c_-|_{T_{k-1}^-} \rho_-|_{T_{k-1}^-} T_k^- -  \Delta t \nabla \cdot \left(\kappa_-|_{T_{k-1}^-} \nabla T_k^- \right) &= \Delta t Q + c_-|_{T_{k-1}^-} \rho_-|_{T_{k-1}^-} T_{n-1}^- \qquad \text { in } \Omega_-,  \\
		T_{k}^- &= T_{k-1}^+ \qquad \text{ on } \gamma.
\end{split}\end{align}
\newline \textbf{Step k+1:} Obtain global solution $T_k^+$:
\begin{align}\begin{split}\label{eq:paralell_update}
		c_-|_{T_{k-1}^+} \rho_-|_{T_{k-1}^+} T_k^+ -  \Delta t \nabla \cdot \left(\kappa_-|_{T_{k-1}^+} \nabla T_k^+ \right) &= \Delta t Q + c_-|_{T_{k-1}^+} \rho_-|_{T_{k-1}^+} T_{n-1}^+ \\ &+ \Delta t \left(\kappa_+|_{T_{k-1}^+} - \kappa_- |_{T_{k-1}^-}\right)\frac{\partial T_{k-1}^-}{\partial \boldsymbol{n}} \xi _{\gamma}   \qquad \text { in } \Omega_+,  \\
		T_{k}^+ &= T_{0} \qquad \text{ on } \Gamma_D.
\end{split}\end{align}
\newline \textbf{Step k+2:} Relaxation: $T_k^+ = \theta T_{k}^+ + (1-\theta)T_{k-1}^+ $, $0<\theta \leq 1$.
\newline \textbf{Step k+3:} Check convergence. If reached, terminate. Else, return to step $k$.
\vspace{5mm}
\\
In \cref{eq:seq_update,eq:paralell_update} $\Gamma_D$ is the Dirichlet boundary which in our case corresponds to the lower surface of the base plate and $T_0=T_{bp}$. In the following, the relaxation parameter $\theta$ is set equal to 1, i.e., no relaxation is applied. 
\par The sequential and parallel formulations are similar, with a key distinction: in step $k+1$, the sequential method utilizes the output $T_k^+$ of step $k$ on the right hand side of \cref{eq:seq_update}, while the parallel method uses the value $T_{k-1}^-$ computed at the previous iteration as force term in \cref{eq:paralell_update}.
While this means that the parallel method may require more iterations to converge, it also means that steps $k$ and $k+1$ are independent of each other, and hence may be computed in parallel. Depending on the computational resources available, the greater parallelism offered by this approach may provide faster computation when compared to the fewer necessary iterations when using the sequential approach. We note that these different algorithms can be seen as Gauss-Seidel (sequential) or Jacobi (parallel) versions of domain decomposition schemes \cite{QV1999, toselli2004domain, DJN2015}. In the following, we adopt the sequential implementation of the two-level method. In fact, due to the limited number of computational cores on our machine, we prefer to exploit the faster convergence rate of the former solution, but for different architectures, the latter might lead to a remarkable computational speed-up.

\subsection{Agglomerated layers activation}\label{ssec:layerActivation}
\begin{figure*}[h]
	\centering
	\includegraphics[width=\textwidth]{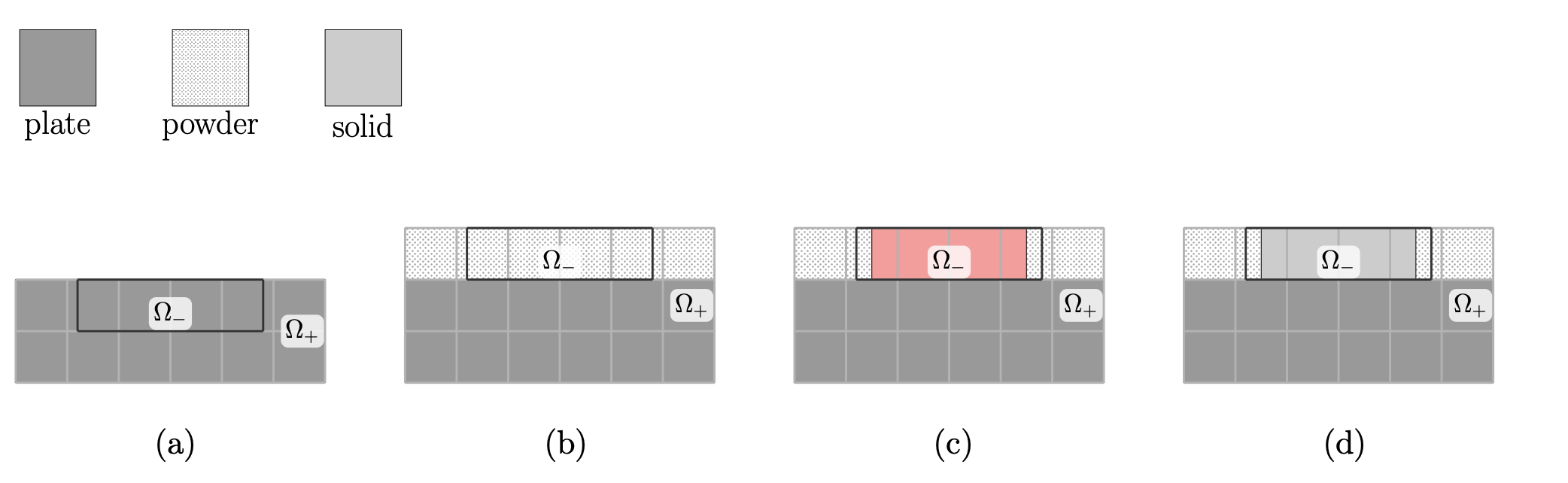}
	\caption{Two level activation strategy. At each activation step a set of 20 agglomerated powder layers is activated. Starting from the build plate at the beginning of the analysis (a), a new agglomerated layer is generated during the layer diffusion step (b), then, during the heating step (c) an equivalent heat source is applied uniformly within the part geometry in the agglomerated layer, and, lastly, the solidified domain is obtained during the cooling step (d). Such a thermal cycle is repeated until the final part is completed.}
	\label{fig:growingDomain}	
\end{figure*}

\cref{fig:growingDomain} depicts the layer-by-layer activation scheme adopted in our implementation. To compute an entire LPBF process at part-scale, we employ an agglomerated layer activation strategy (see, e.g., \cite{Hodge2014,WILLIAMS2018416,carraturo2020modeling}). 
As in \cref{fig:growingDomain}(a), we start first solving a steady-state thermal problem on the build-plate in order to define our initial conditions before the first layer is activated. 
In the present implementation, at each layer activation step, a set of 20 physical layers is inserted into the problem, adding a new layer of global domain elements into the system. The local domain discretization is then shifted upward by the agglomerated layer thickness. Once this mesh update is performed, we compute a layer diffusion step (\cref{fig:growingDomain}(b)), where the residual heat of the previous solution diffuses into the newly activated agglomerated layer. Successively, a heating step occurs, where an equivalent heat source is applied uniformly at the Gauss points of the new layers lying within the geometry of the part. Such a set of Gauss point is evaluated by means of an inside-outside test performed directly on the original .stl file of the part as acquired from a CAD environment. 

In the present work, we define the equivalent heat source $Q$ as:
\begin{equation}
	Q=\dfrac{4 \eta P}{\pi l_d^2 t_a}
	\label{eq:heatSource}
\end{equation}
with $\eta$ the absorptivity of the material, $P$ the nominal power of the laser beam,  $l_d$ the laser spot size, and $t_a$ the thickness of the activated layer.
This definition corresponds to a local, high-intensity heat source instantly activated over the entire portion of the printed geometry in the new agglomerated layer, as depicted in \cref{fig:growingDomain}(c). Therefore, in our implementation the heating activation step lasts few microseconds due to the high energy input resulting from the adopted physical model. In the following, we evaluate the heating activation time step increment $\Delta t_A$ as follows:
\begin{equation}
	\Delta t_A = \dfrac{2l_d}{v}
	\label{eq:HET}
\end{equation}
with $v$ the laser velocity.

The last step of our thermal cycle is the cooling step, wherein the heat source is set equal to zero and the heat is dissipated by conduction through the solidified domain and the surrounding powder and by convection and radiation through the upper surface as described in \cref{ssec:thermalBC}. In the following numerical analysis, a backward Euler implicit time integration scheme is adopted to integrate over time \cref{eq:thermalEquation}.

We employ the so-called \textit{birth-death} element activation scheme (\citet{MARTUKANITZ2014}). In fact, new elements in the global domain are generated as each agglomerated layer is activated. The \textit{birth-death} element activation scheme is chosen since in a structured, cartesian grid-like discretization - such as the one employed by the two-level method - the new element layer generation procedure is straightforward \cite{carraturo2020modeling}. About half of the computational resources are thereby saved during the course of a complete analysis, since no space has to be previously allocated for inactive elements, as required, e.g., by the quiet element method~\citep{Michaleris2014}.

\subsection{Material Properties}\label{ssec:materialProperties}
\begin{center}
	\begin{table*}
		\csvreader[head to column names, 
		tabular=llll, 
		table head=\toprule{Temperature [$^{\circ}$C]} & {Conductivity $k_{sol}$ [W/m/K]} & { Heat capacity $c_{sol}$ [J/kg/K]} & {Density $\rho_{sol}$ [kg/m$^3$]}\\\midrule, 
		table foot=\bottomrule]
		{plots/data/SS316L.csv}{}%
		{\T & \TC & \SHC & \D}
		\caption{Temperature-dependent material properties of solid SS 316L from~\citep{truman2009net}. Values at temperatures higher than the maximum temperature (1400$^{\circ}$C) have been extrapolated. \label{table:SS316L}}
	\end{table*}
\end{center}
In the following, we consider components printed with stainless steel 316L (SS 316L). Such a metal alloy is one of the most widespread materials for metal AM products, thus its material properties are well-known and readily available, even for the high temperatures involved in LPBF processes. In~\cref{table:SS316L} we report the temperature-dependent material properties for solid SS 316L as measured by~\citet{truman2009net}. Since we include powder in the computational domain, we need to define thermal properties for the powder domain region as well. The relation between powder and solid density and heat capacity is then defined as follows:
\begin{align}
	\rho_{pow} &=(1-\phi)\rho_{sol}, \\
	c_{pow} &= c_{sol},
\end{align}
whereas, assuming spherical particles, the relationship between powder and solid thermal conductivity is given by~\citet{sih2004prediction}:
\begin{equation}
	\begin{split}
	\dfrac{\kappa_{pow}}{\kappa_{gas}}= 
	\left(1 - \sqrt{1-\phi}\right)\left(1+\dfrac{\dfrac{4}{3}   \sigma_{SB}  T^3  d_{pow}}{\kappa_{gas}}\right)
	+\sqrt{1-\phi} \dfrac{2}{1-\dfrac{\kappa_{gas}}{\kappa_{sol}}}	\left(\dfrac{2}{1-\dfrac{\kappa_{gas}}{\kappa_{sol}}} \ln
	\left(\dfrac{\kappa_{sol}}{\kappa_{gas}}\right)-1\right)\\ 
	+\sqrt{1-\phi}\dfrac{\dfrac{4}{3}  \sigma_{SB}  T^3  d_{pow}}{\kappa_{gas}},
	\end{split}
	\label{eq:powderConductivity}
\end{equation}
where $\kappa_{gas}$ is the conductivity of the inert gas (Argon) present in the build chamber, $d_{pow}$ is the average diameter of the powder particles. The adopted values of the above quantities are described in \cref{sec:ResultsAndDiscussion}. 
Due to the immersed nature of the two-level method, in both local and global domain there are elements cut by the boundaries of the solidified domain. 

The immersed boundary method can efficiently deal with discontinuous integration at void-material interfaces \cite{Parvizian2007,Duster2017Ency}, whereas numerical instabilities occur when we need to integrate across an interface between two different materials, such as the powder-solid interface in the thermal model described in \cref{sec:governingEquations}. To overcome such a numerical issue, in \cref{eq:powderConductivity} we introduce a smoothing factor $S$ multiplying the left-hand-side term and defined as follows:
\begin{equation}
	S=(1-\sqrt{\nabla T^T \nabla T}h_{el}\delta),
	\label{eq:smoothingFactor} 
\end{equation}
 where $h_{el}$ is the nominal element length and $\delta$ a scaling factor which we set equal to 0.2. As can be noticed, the smoothing factor $S$ defined in \cref{eq:smoothingFactor} depends on the temperature gradient, and thus it almost vanishes when temperature gradients are close to zero as occurs in the layer diffusion and cooling steps. In contrast, it plays an important role in stabilizing the solution during the heating step, when high temperature gradients are involved.
 
\subsection{Domain discretization}
The core idea behind the two-level method is to separate the local and the global scales using two different domains, both of which are discretized separately using potentially very different element length. Moreover, the immersed nature of this numerical approach allow us to employ simple cartesian grids to discretize both the local and the global domain. 
These general concepts of the method can be also applied in the context of part-scale AM thermal process analyses. In particular, in the article at hand, we employ tetrahedral linear elements with a nominal element length $h_{el}$ of $4\, mm$ for the global domain and $1\, mm$ for the local domain.  

\subsection{Error metrics}\label{ssec:ErrMeas}
To quantify the accuracy of our simulated results with respect to the experimental measurements, we employ two metrics, namely the maximum relative percentage error on the dwell temperature $\epsilon_{max}$ and the sample Pearson correlation coefficient $P$ defined as follows:
\begin{equation}
	\epsilon_{max}\% = \maximum{\dfrac{T_{n,i}-T_{m,i}}{T_{m,i}}}\times 100\% \qquad i=1,..,N_{agg}
	\label{eq:relativeError}
\end{equation}
\begin{equation}
\text{P}(N,M)\% = \dfrac{\sum_{i=1}^{N_{agg}}(T_{n,i}-\bar{T}_{n})(T_{m,i}-\bar{T}_{m})}{\sqrt{\sum_{i=1}^{N_{agg}}(T_{n,i}-\bar{T}_{n})^2\sum_{i=1}^{N_{agg}}(T_{m,i}-\bar{T}_{m})^2}}\times 100\%
\label{eq:correlationCoefficient}
\end{equation}
where  $T_{m,i}$ and $T_{n,i}$ are the dwell temperature values of experimental and numerical results at the $i^{th}$ agglomerated layer, $N_{agg}$ is the number of agglomerated layers, $M$ and $N$ are the set of measured and simulated dwell temperatures, and $\bar{T}_{m}$ and $\bar{T}_{n}$ are the corresponding mean values. Therefore, a correlation $P=100\%$ indicates a perfect correlation between the numerical results and the experimental measurements.

\section{Results and discussion}\label{sec:ResultsAndDiscussion}
All the numerical analyses discussed in this section have been obtained using an \textit{in-house} developed code written in \texttt{FreeFEM++}\cite{MR3043640} on a desktop computer provided with Intel\textsuperscript \textregistered ~Xeon\textsuperscript \texttrademark ~W-2123, CPU@3.9GHz, RAM 256Gb. 

\subsection{Model validation on a 3D printed cylinder}\label{ssec:modelCalibration}
\begin{figure*}[h!]
	\centering
	\subfloat[Stl model of the cylinder used in \cite{williams2019situ} (measures are in millimeters).\label{subfig:CylGeom}]
	{
		\includegraphics[width=0.15\textwidth]{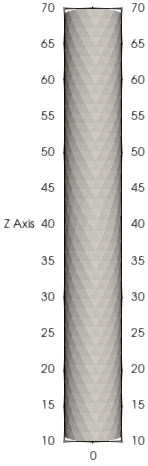}
	}
\hspace{180pt}
	\subfloat[Dwell temperature distribution of the last layer of the cylinder.\label{subfig:CylTemp}]
	{
		\includegraphics[width=0.3\textwidth]{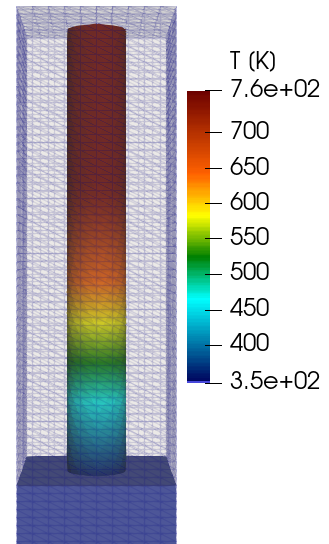}
	}        
	\caption{3D virtual model and thermal simulation results of a cylinder.  \label{fig:Cylinder}}
\end{figure*}

To calibrate and validate the accuracy of the proposed numerical method, we consider the LPBF printing process of a cylinder as reported by Williams \textit{et al.} \cite{williams2019situ} and depicted in \cref{subfig:CylGeom}. Williams and co-workers measure the average temperature on the upper powder layer before laser activation, which they refer to as dwell temperature. According to their experimental findings, such a dwell temperature is controlled by the interlayer cooling time (ILCT) and plays a crucial role in defining the microstructure and thus the mechanical properties of AM parts. 

\begin{figure*}[h!]
\centering
\includegraphics[width=\textwidth]{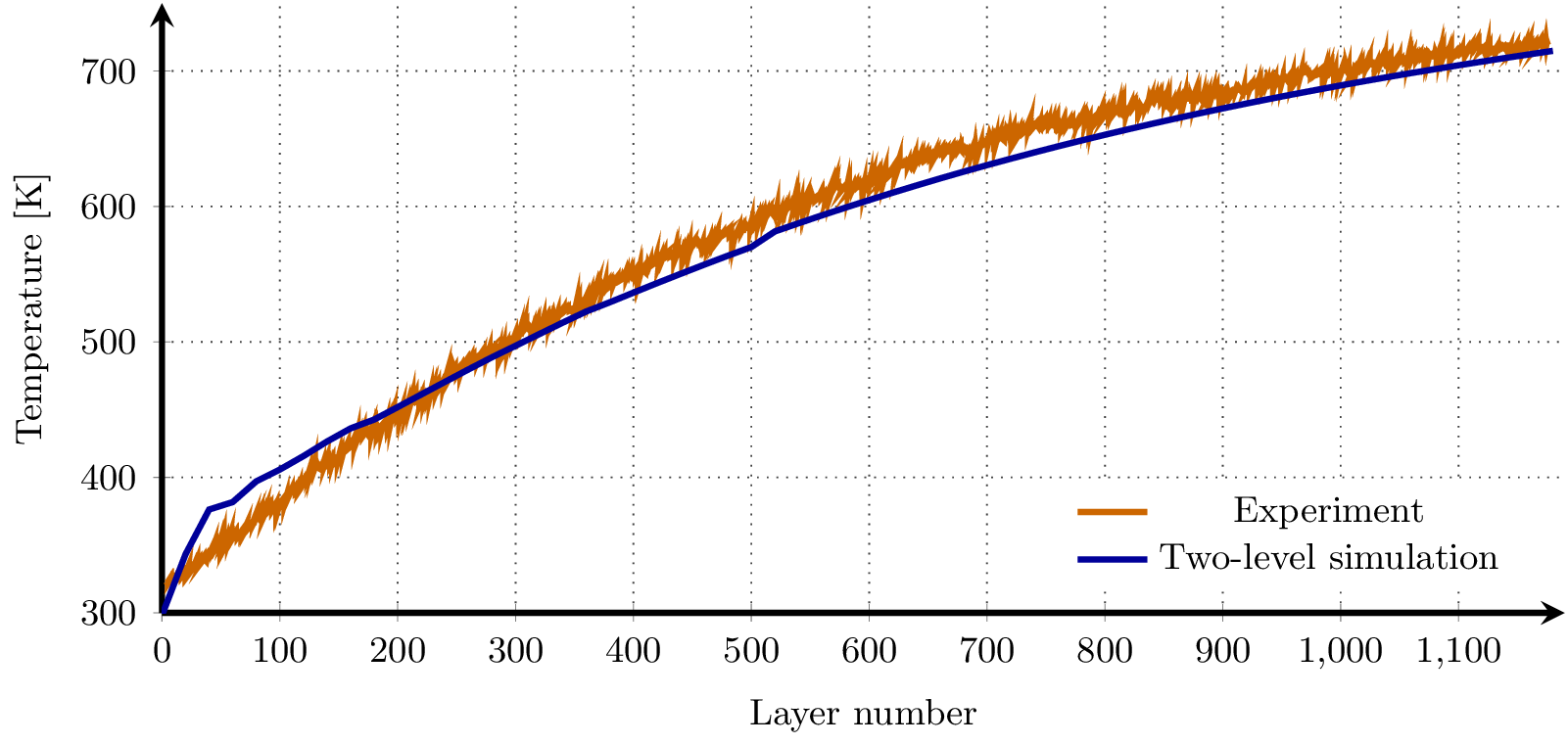}
	\caption{Simulated and measured surface temperature of the cylinder taken just before the laser passes over each layer considering an ILCT of 11s.}
	\label{fig:ILCT_temp}	
\end{figure*}

We have replicated their experiment using the previously presented part-scale two-level method. The simulated dwell temperature evolution at the last layer of the LPBF process is shown in \cref{subfig:CylTemp} using the model parameters reported in \cref{table:calibratedParam}, which have been calibrated following a procedure similar to the one described in \cite{kollmannsberger2019}. 
In \cref{fig:ILCT_temp} we report the simulated and the measured dwell temperatures during the entire 3D printing process. The dwell temperature results obtained by means of the two-level part-scale model shows good agreement with the experimental measurements reported in \cite{williams2019situ}. The maximum relative error $\epsilon_{max}$ is equal to $13.95\%$ and occors at the second agglomerated layer, i.e., in a transitory step, whilst the relative error is kept always below $2.5\%$ after the $15^{th}$ activation step. The Pearson correlation coefficient as defined in \cref{ssec:ErrMeas} is equal to 99.74$\%$, indicating an almost perfect correlation among the measured and the predicted dwell temperature results, making us confident on the reliability of the proposed part-scale thermal model.  

\begin{table}[h!]
	\begin{tabular}{lll}
		\hline
		Parameters [unit] & Values  & Sources \\
		\hline 
		laser power [W]	& 200 & \cite{williams2019situ} \\
		laser spot size [$\mu$m] & 65 & " \\
		laser velocity [mm/s] & 800 & " \\
		powder layer thickness [$\mu$m] & 50 & " \\
		build chamber temperature [$^{\circ}$C] & 25 & " \\
		build plate temperature [$^{\circ}$C] & 80 & " \\
		HET [$\mu$s] & 80 & \refeq{eq:HET} \\
		gas conductivity $\kappa_{gas}$ [W/m/K] & 0.0172 & Calibrated \\
		porosity $\phi$ [-] & 0.35 & " \\
		$h_{pow}$ [W/m/K] & 25 & " \\ 
		$h_{conv}$ [W/m/K] & 0.1 & " \\ 
		emissivity $\varepsilon$ [-] & 0.25 & " \\ 
		absorptivity $\eta$ [-] & 0.7 & " \\ 
		\hline
	\end{tabular}
	\caption{LPBF process and model parameters. \label{table:calibratedParam}}
\end{table}

\subsection{Thermal analysis of a topologically optimized structure}\label{ssec:numericalExample}
\begin{figure}[h!]
	\centering
	\includegraphics[width=\textwidth]{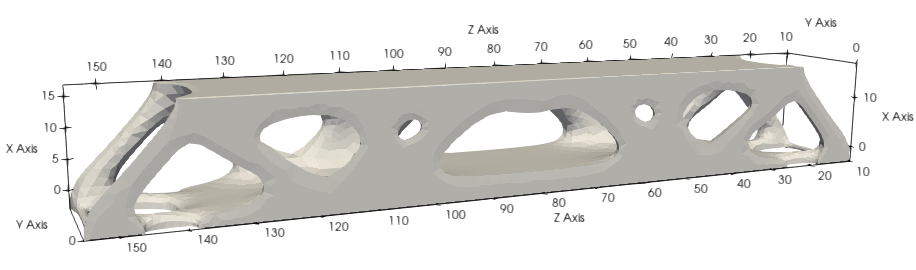}
	\caption{3D virtual model of topologically optimized beam structure taken from \cite{Carraturo2021GAMM} (measures are in millimeters).}
	\label{fig:TOD_solid_full}	
\end{figure}

\begin{figure*}[h!]
	\centering
	\subfloat[Agglomerated layer 52.\label{subfig:TOT_52}]
	{
		\includegraphics[width=0.65\textwidth]{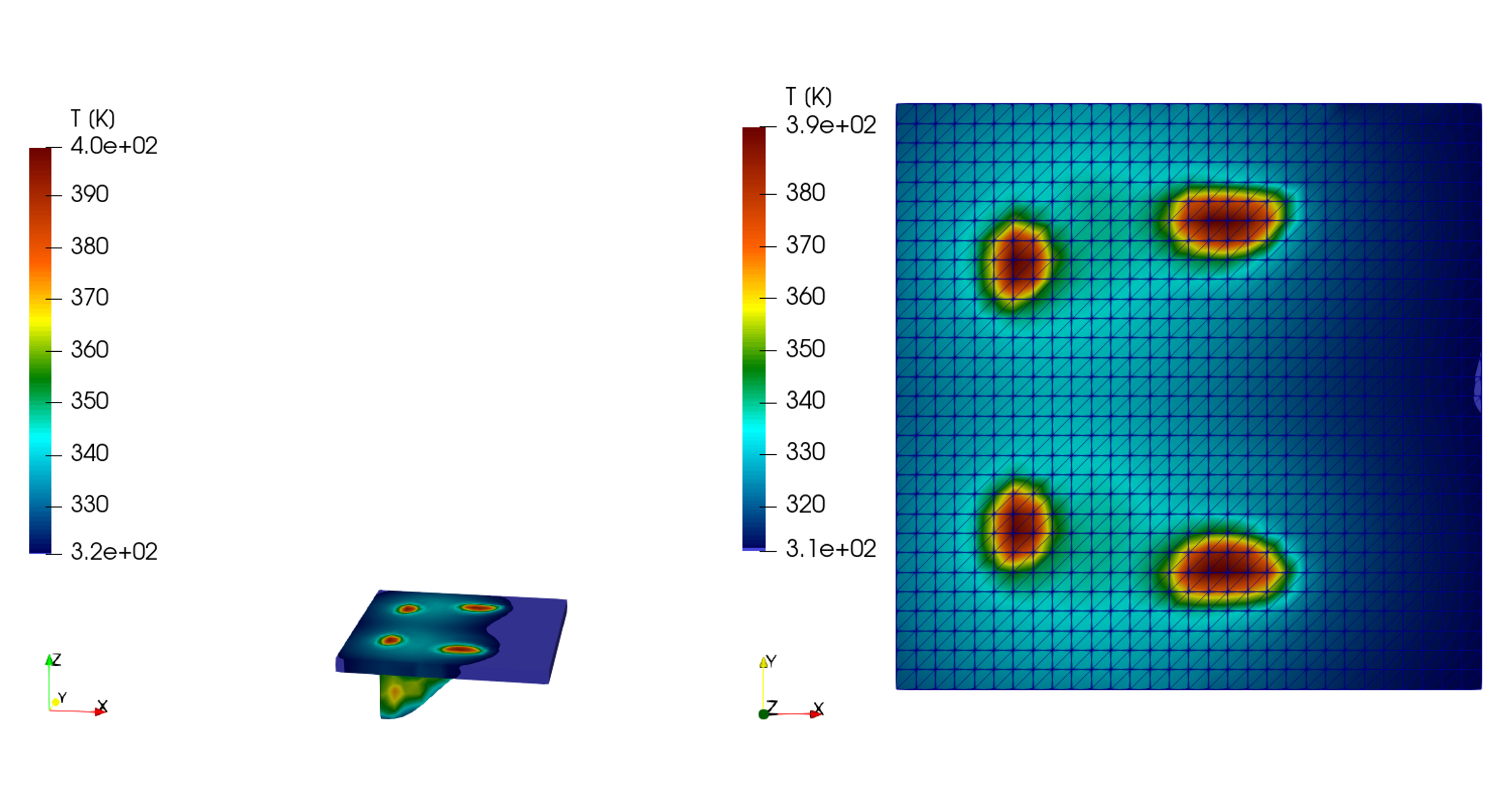}
	}
\\
	\subfloat[Agglomerated layer 307.\label{subfig:TOT_307}]
	{
		\includegraphics[width=0.65\textwidth]{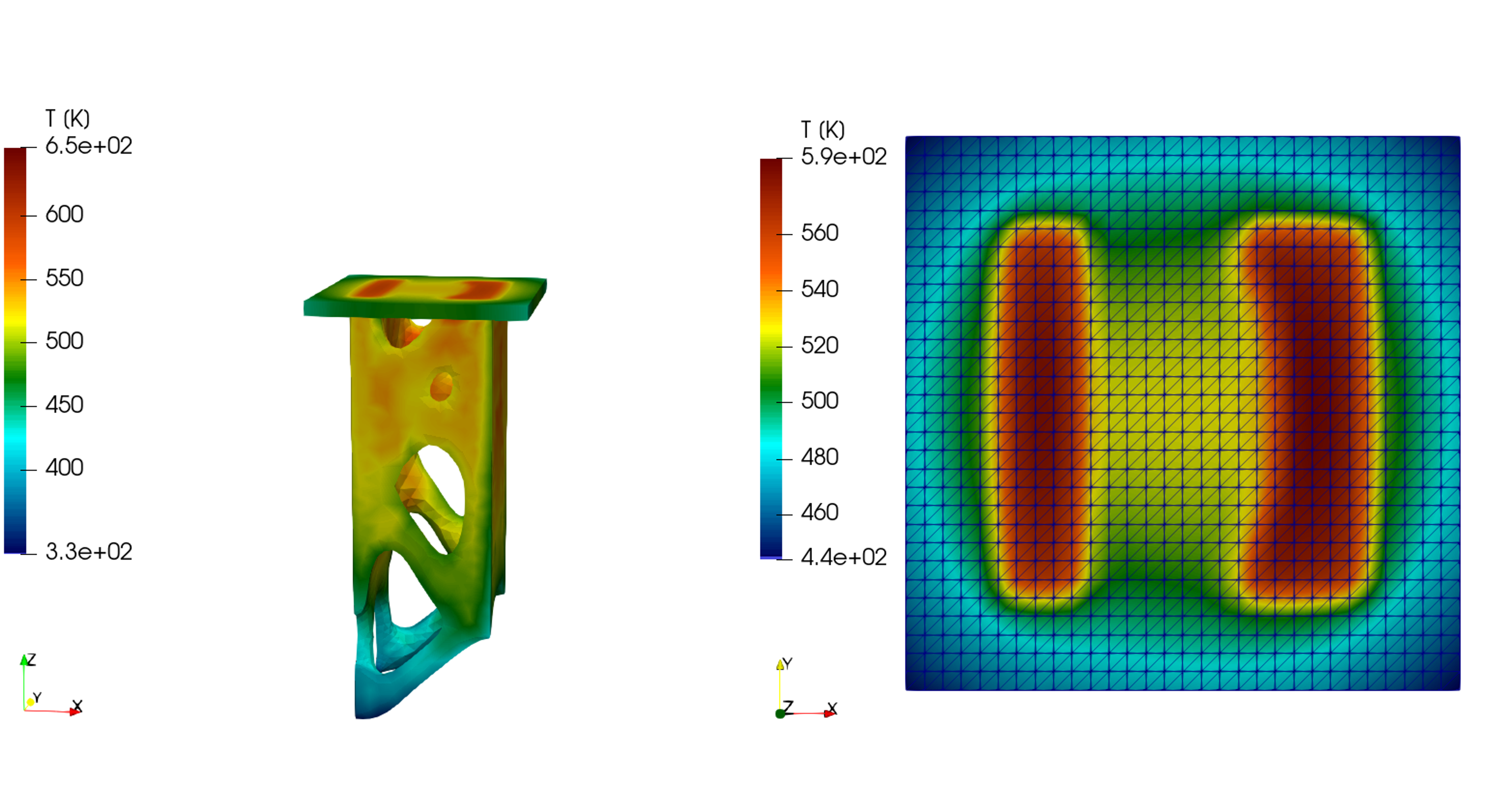}
	}  
\\
\subfloat[Agglomerated layer 453.\label{subfig:TOT_453}]
{
	\includegraphics[width=0.65\textwidth]{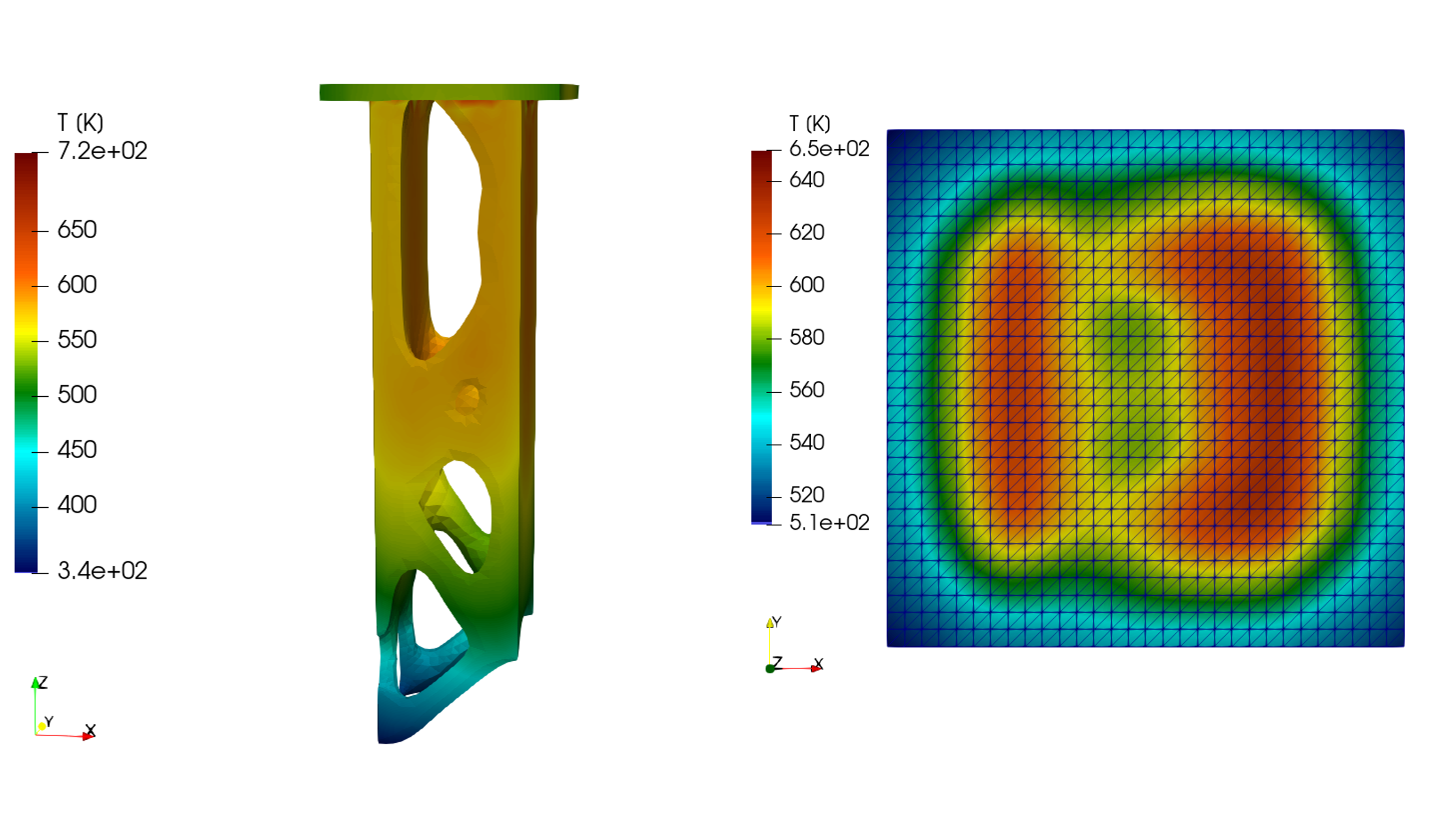}
}       
	\caption{Temperature distribution at different layers. On the right it is reproduced the simulated temperature field on both the growing printed geometry and the local domain; on the left the top view showing the local domain temperature distribution and the corresponding domain discretization.  \label{fig:TemperatureOptimizedMBB}}
\end{figure*}

In order to test the capability of the proposed numerical framework to deal with complex geometries, we compute the part-scale thermal analysis of a topologically optimized beam structure. The optimized structure depicted in \cref{fig:TOD_solid_full} was obtained in \cite{Carraturo2021GAMM} using a phase-field-based topology optimization procedure as presented in \cite{Auricchio2020,Carraturo2019TopOpt}.
In \cref{fig:TemperatureOptimizedMBB} we report the temperature field distribution on the topologically optimized part and on the local domain at three different steps of the process as well as the corresponding local domain mesh. Observing the results reported in \cref{fig:TemperatureOptimizedMBB}, we can observe the two-level domain concept in action. On the one hand, the local domain allows to capture the small geometrical details of the optimized structure with a very high resolution, on the other hand, the coupled global domain solution approximates the temperature field in the remaining regions of the part, substantially accounting for the thermal history at the scale of the entire component. Such an approach does not require any meshing of the considered geometry, since the domain is treated in an implicit fashion, making the analysis setup straightforward and suitable to deal with geometric flaws and complex geometries.  

\section{Conclusions}\label{sec:conclusions}

In the present contribution, we have applied the two-level framework, first introduced in \cite{viguerie2020fat, viguerie2021numerical} for high-fidelity thermal models, to the analysis at the scale of an entire 3D printed component. Such an extension of the physical model in the implemented two-level framework allows us to simulate the complete printing process of parts with complex geometries, without requiring any conformal mesh generation step, leading to a smooth workflow from the CAD model to the thermal analysis of AM components. The presented numerical results show the effectiveness of the proposed approach.
\par There are several potential perspectives for the present work. First, the structured discretization grids of both the local and the global domain call for a massive code parallelization. Such a code optimization step would lead to a remarkable speed-up in the overall computational time, fully exploiting the potentiality of the proposed numerical framework and opening the possibility to compute thermal analysis of an entire build plate. Secondly, we aim at extending the proposed part-scale thermal model to a multi-physics setting, considering a weakly-coupled thermo-mechanical problem formulation. 
\par As shown in our numerical examples, there is ample potential to use the machinery shown here to better unify the development and analysis of topologically optimized components, as the potential reduction in cost afforded by the two-level method can greatly accelerate forward solves during the optimization process.
In fact, the optimization routine may require the numerical solution of many thousands of problems. In this aim, the extension of the machinery shown here to include data-driven and reduced order methods, including POD- and DMD- based methods, as well as concepts from scientific machine learning, may allow for even more efficient numerical solutions. 

\section*{Acknowledgments}
This work was partially supported by the Italian Minister of University and Research through the MIUR-PRIN projects "A BRIDGE TO THE FUTURE” (No. 2017L7X3CS) and "XFAST-SIMS" (no. 20173C478N).
%% References -------------------------------------
%\bibliographystyle{apalike}
\bibliographystyle{elsarticle-num-names}
\biboptions{sort&compress}
\bibliography{TwoLevel}

\end{document}